\documentclass{amsart}
\usepackage{amsmath,amsthm,amssymb,amsfonts,graphicx,mathrsfs,bm,amscd,latexsym}
\usepackage{color,enumitem,hyperref}
\usepackage{blindtext}
\usepackage{textcomp}
\usepackage[all]{xy}

\newtheorem{maintheorem}{Theorem}
\newtheorem{maincor}{Corollary}

\newtheorem{theorem}{Theorem}[section]
\newtheorem{proposition}[theorem]{Proposition}
\newtheorem{lemma}[theorem]{Lemma}

\def\R{\mathbb{R}}

\def\N{\mathbb{N}}

\def\S{\mathbb{S}}
\def\W{\mathcal{W}}

\def\g{\gamma}

\def\D{\Delta}

\def\C{\mathcal{C}}

\DeclareMathOperator{\diam}{diam}
\DeclareMathOperator{\dist}{dist}

\numberwithin{equation}{section}

\begin{document}

\title{Uniformization of Cantor sets with bounded geometry}

\address{Department of Mathematics, The University of Tennessee, Knoxville, TN  37966}

\email{vvellis@utk.edu}

\thanks{The author was partially supported by the Academy of Finland project 257482 and NSF DMS grant 1952510.}
\author{Vyron Vellis}

\date{\today}


\begin{abstract}
In this note we provide a quasisymmetric taming of uniformly perfect and uniformly disconnected sets that generalizes a result of MacManus \cite{MM2} from 2 to higher dimensions. In particular, we show that a compact subset of $\R^n$ is uniformly perfect and uniformly disconnected if and only if it is ambiently quasiconformal to the standard Cantor set $\mathcal{C}$ in $\R^{n+1}$.
\end{abstract}

\maketitle

\section{Introduction}\label{sec:cantor}

The \emph{(quasisymmetric) uniformization problem} asks for necessary and sufficient conditions under which a metric space $X$ is quasisymmetrically homeomorphic to a ``standard" space $X_0$. Roughly speaking, quasisymmetric homeomorphisms are a generalization of conformal maps which preserve relative distances; see Section \ref{sec:prelim} for precise definitions. The uniformization problem has been extensively studied in literature for a variety of ``standard" spaces: the unit circle $\mathbb{S}^1$ \cite{TuVa}, the unit sphere $\mathbb{S}^2$ \cite{BonK}, geodesic trees \cite{BM}. See also \cite{Bonk6} for a general overview.

Perhaps the most simple metric spaces from a topological point of view are Cantor sets, i.e., homeomorphic images of the standard ternary Cantor set $\mathcal{C}$. Brouwer's topological characterization of Cantor sets \cite[Theorem 7.4]{Kechris} states that a metric space is a Cantor set if and only if it is compact, perfect, and totally disconnected. 

In \cite{DSbook}, David and Semmes solved the uniformization problem in the case that the standard space $X_0$ is $\mathcal{C}$. Contrary to Brouwer's uniformization, in the quasisymmetric case one has to assume some quantitative versions of perfectness and total disconnectedness. 

A closed non-degenerate metric space $X$ is called \emph{uniformly perfect} if there exists a constant $c\geq 1$ such that for all $x \in X$ and all $r \in (0,\diam{X})$, there exists a point in $B(x,r) \setminus B(x,r/c)$. Every uniformly perfect space is perfect; on the other hand, the planar set $[0,1]\times\bigcup_{n=1}^{\infty}\{n!\}$ is perfect but not uniformly perfect.  

A metric space $X$ is \emph{uniformly disconnected} if there is a constant $c\geq 1$ such that for all $x\in X$ and all positive $r<\frac{1}{4}\diam{X}$, there exists $E\subset X$ containing $x$ such that $\diam{E} \leq r$ and $\dist(E,X\setminus E) \geq r/c$. All uniformly disconnected spaces are totally disconnected; on the other hand, $\mathbb{Z}$ is totally disconnected but not uniformly disconnected.

Based on these scale-invariant notions, David and Semmes uniformization is as follows.

\begin{theorem}[{\cite[Proposition 15.11]{DSbook}}]\label{thm:UDQS}
A metric space is quasisymmetrically homeomorphic to $\C$ if and only if it is compact, doubling, uniformly disconnected and uniformly perfect.
\end{theorem}

Recall that a metric space $X$ is doubling if there exists a constant $C>1$ such that for all $x\in X$ and all $r>0$, the ball $B(x,r)$ can be covered by at most $C$ many balls of radius $r/2$. Since all Euclidean spaces $\R^n$ are doubling, the doubling condition in Theorem \ref{thm:UDQS} can be dropped if $X \subset \R^n$. 

Later, MacManus \cite{MM2} proved a stronger uniformization result for Cantor sets contained in $\R^2$. The improvement here is that the quasisymmetric homeomorphism can be in fact assummed to be defined on the ambient space $\R^2$ and not just $\C$.

\begin{theorem}[{\cite[Theorem 3]{MM2}}]\label{thm:MM}
For a compact set $X\subset \R^2$ there exists a quasisymmetric mapping $F:\R^2 \to \R^2$ with $F(\C) = X$ if and only if $X$ is uniformly perfect and uniformly disconnected. 
\end{theorem}

Note that Theorem \ref{thm:MM} is false in $\R^3$ and in $\R^4$ due to the existence of a self-similar wild Cantor set in $\R^3$ called \emph{Antoine necklace}; see \cite[pp. 70--75]{Daverman}. By self-similarity, this set is both uniformly perfect and uniformly disconnected, but there exists no homeomorphism of $\R^3$ (let alone a quasisymmetric homeomorphism) that maps this set onto $\C$. See also  \cite[Appendix A]{PW2} for recent examples in $\R^4$.

In our main result, we provide a quasisymmetric taming of Cantor sets with bounded geometry. In particular, we show that by increasing the dimension by $1$, MacManus' result generalizes to all dimensions $n\geq 3$.

\begin{maintheorem}\label{thm:main}
Let $n\in\N$. 
For a compact set $X\subset \R^n$ there exists a quasisymmetric map $F:\R^{n+1} \to \R^{n+1}$ with $F(\C) =X$ if and only if $X$ is uniformly perfect and uniformly disconnected.
\end{maintheorem}

Here and for the rest, given $n\in\N$, we identify $\R^n$ with the plane $\R^n \times\{0\} = \{(x_{1},\dots,x_n,0) : x_i\in\R\} \subset \R^{n+1}$. 

One application of Theorem \ref{thm:main} is in conformal dynamics. A \emph{uniformly quasiregular map} (abbv. \emph{UQR} map) $f:\R^n \to \R^n$ is a map in the Sobolev space $W^{1,n}_{\text{loc}}(\R^n)$ for which there exists $K\geq 1$ such that for any $m\in\N$, the $m$-th iterate $f_m = f\circ\cdots\circ f$ satisfies 
\[ |f_m'(x)| \leq K J_{f_m}, \qquad \text{for a.e. $x\in\R^n$}.\]
Due to the rigidity of conformal maps in dimensions $n\geq 3$, UQR maps play the role of holomorphic maps in the study of conformal dynamics in higher dimensions. A well known problem in conformal dynamics is the characterization of closed sets in $\R^n$ that arise as Julia sets of UQR maps. Iwaniecz and Martin  \cite{IM} showed that the Cantor set $\mathcal{C}$ is a Julia set of a UQR map of $\R^2$, and Fletcher and Wu \cite{FlWu} showed that the Antoine necklace is a Julia set of a UQR map of $\R^3$. Note that both $\mathcal{C}$ and the Antoine necklace are uniformly perfect and uniformly disconnected because they are self-similar. In general, it is unknown whether self-similar Cantor sets in dimensions $n\geq 3$ are always Julia sets of UQR maps. In \cite{FV} we apply Theorem \ref{thm:main} to show that every uniformly perfect and uniformly disconnected subset of $\R^n$, $n\geq 3$, is the Julia set of a UQR map of $\R^{n+1}$.

Moreover, Theorem \ref{thm:MM} and Theorem \ref{thm:main} yield the following quasiconformal embedding result for uniformly disconnected sets.

\begin{maincor}\label{cor:2}
Let $n\geq 2$ be an integer and let $X\subset \R^{n}$ be a bounded uniformly disconnected set. There exists a quasisymmentric homeomorphism of $\R^{N}$ that maps $X$ into $\mathcal{C}$, where $N=2$ if $n=2$ , and $N=n+1$ if $n\geq 3$.
\end{maincor}

Corollary \ref{cor:2} has an application in hyperbolic geometry. If $X\subset \S^2$ is a Cantor set, then by the Uniformization Theorem, $S:= \S^2 \setminus X$ is necessarily a hyperbolic Riemann surface. Hence, $S$ has a pants decomposition, that is, $S = \bigcup_{i=1}^{\infty} P_i$, where each $P_i$ is a topological sphere with three disks removed. The collection of boundary curves of the pairs of pants, called the \emph{cuffs of the decomposition}, may be enumerated by $(\alpha_j)_{j=1}^{\infty}$. Each $\alpha_j$ is a simple closed curve on $S$ and generates a class $[\alpha_j]$ of simple closed curves that are freely homotopic to $\alpha_j$. Denote by $\ell_S[\alpha_j]$ the infimum of hyperbolic lengths of curves in $[\alpha_j]$. Sugawa  \cite{Sugawa} proved that the Cantor set $X$ is uniformly perfect if and only if $\inf \ell[\alpha_j] >0$. In an upcoming paper with Fletcher, we apply Corollary \ref{cor:2} to show that a similar statement holds for uniformly disconnected sets: a Cantor set $X\subset \R^2$ is uniformly disconnected if and only if $\sup \ell[\alpha_j] < \infty$.

By properties of quasisymmetric homeomorphisms, one direction of Theorem \ref{thm:main} is clear. Namely, if there exists a quasisymmetric map $F:\R^{n+1} \to \R^{n+1}$ with $F(\C) =X$, then $X$ is compact, uniformly perfect and uniformly disconnected. For the converse, which is the content of this paper, we use the existence of a quasisymmetric homeomorphism $f:\mathcal{C} \to X$, and we extend this mapping quasisymmetrically to $\R^{n+1}$. In \textsection\ref{sec:prelim} we give some basic definitions,  in \textsection\ref{sec:schoenflies} we discuss some simple bi-Lipschitz Schoenflies theorems in higher dimensions, and in \textsection\ref{sec:proof} we prove Theorem \ref{thm:main} and Corollary \ref{cor:2}.

\section{Preliminaries}\label{sec:prelim}

For $n\in\N$, a point $x\in \R^n$, and $r>0$ we denote by $B^n(x,r)$, $\overline{B}^n(x,r)$ the open and closed, respectively, balls centered at $x$ and with radius $r$.

A homeomorphism $f\colon D\to D'$ between two domains in $ \mathbb{R}^n$ is called $K$-\emph{quasiconformal} for some $K\geq 1$ if, for all $x\in D$, $f$ satisfies the distortion inequality
\[\limsup_{r\to 0} \frac{\sup_{y\in\partial B^n(x,r)}|f(x)-f(y)|}{\inf_{y\in\partial B^n(x,r)}|f(x)-f(y)|} \leq K.\]

A homeomorphism $f:(X,d_X) \to (Y,d_Y)$ between metric spaces is said to be $\eta$-\emph{quasisymmetric} if there exists a homeomorphism $\eta \colon [0,\infty) \to [0,\infty)$ such that for all $x,a,b \in X$ with $x\neq b$ 
\[ \frac{d_Y(f(x),f(a))}{d_Y(f(x),f(b))} \leq \eta \left ( \frac{d_X(x,a)}{d_X(x,b)} \right ). \]

A quasisymmetric mapping between two domains in $\R^n$ is quasiconformal. The converse holds true for the smaller class of uniform domains which contains $\R^n$.  For a systematic treatment of quasiconformal mappings see \cite{Vais1}.

It follows easily from the definitions above that quasisymmetric maps preserve the notions of uniform perfectness and uniform disconnectedness quantitatively.

\begin{lemma}
If $f : X \to Y$ is $\eta$-quasisymmetric and $X$ is $c$-uniformly perfect (resp. $c$-uniformly disconnected), then $Y$ is $c'$-uniformly perfect (resp. $c'$-uniformly disconnected) with $c'$ depending only on $\eta$ and $c$.
\end{lemma} 

A map $f\colon X \to Y$ between metric spaces is \emph{$L$-bi-Lipschitz} for some $L \geq 1$ if 
\[ L^{-1}d_X(x,y) \leq d_Y(f(x),f(y)) \leq Ld_X(x,y)\] 
for all $x,y \in X$. Note that an $L$-bi-Lipschitz mapping is $L^2t$-quasisymmetric. 

A weaker notion of bi-Lipschitz mappings is that of \emph{bounded length distortion} (\emph{BLD}) mappings. A mapping $f:(X,d_X) \to (Y,d_Y)$ is $L$-BLD if there exists $L \geq 1$ such that 
\[ L^{-1}\ell(\g) \leq \ell(f(\g)) \leq L\ell(\g)\] 
for all paths $\g : [0,1] \to X$. Here and for the rest, $\ell$ denotes the length of a path. Clearly, $L$-bi-Lipschitz mappings are $L$-BLD mappings but BLD mappings need not be bi-Lipschitz even if they are homeomorphisms. However, BLD homeomorphisms between convex spaces are bi-Lipschitz. 

\begin{lemma}\label{lem:BLD}
Let $f: X \to Y$ be an $L$-BLD homeomorphism between two convex metric spaces. Then $f$ is $L$-bi-Lipschitz.
\end{lemma}

An embedding $f:(X,d_X) \to (Y,d_Y)$ is a \emph{$(\lambda,L)$-quasisimilarity} for some $\lambda>0$ and $L\geq 1$ if 
\[ L^{-1}\lambda d_X(x,y) \leq d_Y(f(x),f(y)) \leq L\lambda d_X(x,y) \qquad\text{for all}\quad x,y \in X.\] 
Note that $(\lambda,1)$-quasisimilarities are similarities with scaling factor $\lambda$, while $(1,L)$-quasisimilarities are $L$-bi-Lipschitz, and $(1,1)$-quasisimilarites are isometries.

While similarities preserve relative distances between nondegenare sets, quasisymmetric maps quasi-preserve relative distances between nondegenare sets. Specifically, if $f: X \to Y$ is $\eta$-quasisymmetric and $E, E' \subset X$ are nondegenerate closed sets, then
\begin{equation}\label{eq:relQS}
\frac{1}{2}\phi\left ( \frac{\dist(E,E')}{\diam{E}}\right ) \leq \frac{\dist(f(E),f(E'))}{\diam{f(E)}}  \leq \eta\left ( 2\frac{\dist(E,E')}{\diam{E}} \right )
\end{equation}
where $\phi(t) = (\eta(t^{-1}))^{-1}$; see for example \cite[p. 532]{Tyson}. Moreover, if $f: X \to Y$ is $\eta$-quasisymmetric and $A\subset B \subset X$ are such that $0<\diam{A}\leq \diam{B} < \infty$, then $\diam{f(B)}$ is finite and
\begin{equation}\label{eq:relQS2}
\left(2\eta\left(\frac{\diam{B}}{\diam{A}}\right)\right)^{-1} \leq \frac{\diam{f(A)}}{\diam{f(B)}} \leq \eta\left(2\frac{\diam{A}}{\diam{B}}\right).
\end{equation}
For the proof of (\ref{eq:relQS2}) see \cite[Proposition 10.8]{Heinonen}.

\section{Bi-Lipschitz Schoenflies results}\label{sec:schoenflies}

The classical Sch\"onflies theorem states that every embedding of $\S^1$ in $\R^2$ extends to a homeomorphism of $\R^2$. Tukia \cite{TukiaBLExt} proved a bi-Lipschitz version of Sch\"onflies theorem.

\begin{theorem}\cite{TukiaBLExt}\label{thm:BAT}
If $f:\S^1\to\R^2$ is 
$L$-bi-Lipschitz
, then $f$ extends 
$L'$-bi-Lipschitz
to $\R^2$ with 
$L'$ depending only 
on $L$.
\end{theorem}

It is well known that in higher dimensions Theorem \ref{thm:BAT} fails even under strong topological assumptions. In particular, Tukia \cite[\textsection15]{TukiaBLExt} constructed a bi-Lipschitz embedding of $\S^2$ into $\R^3$ that can be extended as a homeomorphism of $\R^3$ but not as a quasisymmetric (let alone bi-Lipschitz) homeomorphism of $\R^3$.

Theorem \ref{thm:BAT} was generalized for annuli in $\R^2$ and annuli in higher dimensions with controlled topology and geometry by Sullivan \cite{Su}. For the proof of the following theorem see Theorem 3.12, Theorem 5.8 and \textsection5.9 in \cite{TukiaVais-LIPandLQCextension}.

\begin{theorem}\label{thm:annulus}
If $D_1\subset D_2 \subset \R^n$ and $D_1'\subset D_2' \subset \R^n$ are $n$-cubes satisfying 
\begin{equation}\label{eq:annulusthm}
C^{-1}\diam{D_2} \leq \diam{D_1} \leq C \dist(\partial D_2, \overline{D_1})
\end{equation}
for some $C>1$ and if $f : \partial D_1 \cup \partial D_2 \to \partial D_1' \cup \partial D_2'$ is $L$-bi-Lipschitz that admits a homeomorphic extension on $\overline{D_2\setminus D_1}$, then $f$ admits an $L'$-bi-Lipschitz extension on $\overline{D_2\setminus D_1}$ with $L'$ depending only on $L$ and $n$.
\end{theorem} 

In this section, we work with a much simpler setting. For $d>1$ and $n\in\{2,3,\dots\}$ denote by $\mathscr{C}_n(d)$ the collection of domains $U\subset \R^n$ whose boundary components are boundaries of $n$-cubes of diameters and mutual distances bounded below by $d^{-1}\diam{U}$. 

The next proposition is the main result of this section.

\begin{proposition}\label{prop:ext-ndim}
Let $U\in\mathscr{C}_n(d)$ and $f:\partial U \to \R^n$ be an $L$-bi-Lipschitz map that is a similarity on each component of $\partial U$ and that extends homeomorphically to $\overline{U}$. Then $f$ extends $L'$-bi-Lipschitz to $\overline{U}$ with $L'$ depending only on $L$, $d$ and $n$.
\end{proposition}

%
 We start with the simple observation that every domain in $\mathscr{C}_n(d)$ has a finite number of boundary components.

\begin{lemma}\label{lem:numberofcomp}
Let $n\in \{2,3,\dots,\}$ and $d>1$. Every domain $U \in\mathscr{C}_n(d)$ contains at most $N$ boundary components with $N$ depending only on $n,d$.
\end{lemma}

\begin{proof}
Let $D_1,\dots,D_m$ be some bounded components of $\R^n \setminus \overline{U}$. For each $i\in\{1,\dots,m\}$, choose $x_i\in D_i$. Then, for each distinct  $i,j$, we have $x_i \in \overline{B}^n(x_1,\diam{U})$ and $|x_i-x_j| \geq d^{-1}\diam{U}$. By the doubling property of $\R^n$, we have that $m\leq N_0d^{-n}$ for some universal $N_0$.
\end{proof} 

Below, for a domain $U \in\mathscr{C}_n(d)$, we write $U = U(D_0;D_1,\dots,D_m)$ if $D_0,\dots,D_m$ are open cubes with 
\begin{enumerate}
\item $U = D_0 \setminus \bigcup_{i=1}^m \overline{D_i}$;
\item $\overline{D_1},\dots,\overline{D_m}$ are contained in $D_0$ and are mutually disjoint;
\item for all $i\in\{1,\dots,m\}$ we have $\diam{D_i} \geq d^{-1}\diam{D_0}$;
\item for all distinct $i,j \in \{0,\dots,m\}$ we have $\dist(\partial D_i, \partial D_j) \geq d^{-1}\diam{D_0}$.
\end{enumerate}

For the rest of Section 3 we denote by $\mathcal{U}_0$ the open $n$-cube $(-1,1)^n$ in $\R^n$. For each $m\in\N$ and $k\in \{1,\dots,m\}$ denote
\[\mathcal{S}_{m,k} := \left [\frac{4k-2m-3}{2m+1},\frac{4k-2m-1}{2m+1}\right ]\times\left[\frac{-1}{2m+1},\frac{1}{2m+1}\right]^{n-1} \subset \R^n\]
and 
\[ \mathcal{U}_m = \mathcal{U}_0 \setminus \bigcup_{k=1}^m \mathcal{S}_{m,k}.\]


In the next lemma, we show that every domain $U \in \mathscr{C}_n(d)$ is quasisimilar to $\mathcal{U}_m$ for some $m\in\N$. This allows us to reduce the proof of Proposition \ref{prop:ext-ndim} to the case $U = \mathcal{U}_m$.

\begin{lemma}\label{lem:standdom}
For $n\in \{2,3,\dots\}$ and $d>0$ there exist $L\geq 1$ depending only on $n,d$ with the following property. If $U = U(D_0;D_1,\dots,\cdots,D_m)$ is in $\mathscr{C}_n(d)$ and has diameter equal to $1$, then there exists an $L$-bi-Lipschitz homeomorphism $f:\overline{U} \to \overline{\mathcal{U}_m}$ with $f(\partial D_0) = \partial \mathcal{U}_0$ and $f(\partial D_i) = \partial \mathcal{S}_{m,i}$ for $i=1,\dots,m$.
\end{lemma} 

For the proof of Lemma \ref{lem:standdom}, we say that two closed cubes $D,D' \subset \R^n$ are \emph{concentric} if there exist $r,r'>0$ and an isometry of $\R^n$ that maps $D$ onto $[-r,r]^n$, and $D'$ onto $[-r',r']^n$. The \emph{center} of a closed cube $D\subset \R^n$ is the unique point $x_0 \in D$ such that there exists $r>0$ and an isometry $f$ of $\R^n$ that maps $D$ onto $[-r,r]^n$ and $x_0$ onto the origin.

Recall also that a \emph{dyadic $n$-cube} $D\subset \R^n$ is a set 
\[ D = [i_12^{k}, (i_1+1)2^k]\times \cdots\times [i_n2^{k}, (i_n+1)2^k]\] 
where $k,i_1,\dots,i_n \in \mathbb{Z}$.




\begin{proof}[{Proof of Lemma \ref{lem:standdom}}]
%
By Lemma \ref{lem:numberofcomp}, $m\leq N$ for some $N$ depending only on $d$ and $n$. Applying an isometry of $\R^n$, we may assume that $D_0 = \mathcal{U}_0$.

Let $k_1\in\N$ be an integer such that
\[ \tfrac12\min\{(4d)^{-1},(4m+2)^{-1}\} \leq 2^{-k_1} < \min\{(4d)^{-1},(4m+2)^{-1}\}.\]
There exists $x_0 \in U$ such that $\dist(x_0,\partial U) \geq (2d)^{-1}$ and by design of $k_1$, there exists a dyadic cube $\D \subset U$ of side length $2^{-k_1}$ and $\dist(\D,\partial U) \geq 2^{-k_1}$. Let $k_2 \in \N$ be an integer such that 
\[ m< 2^{k_2} \leq 2m.\] 
Divide $\D$ into $2^{k_2}$ dyadic cubes of side-length $2^{-k_1-k_2}$ and let $\D_1',\dots,\D_m'$ be $m$ of them. Inside each $\D_i'$ fix a dyadic cube $\D_i \subset \D_i'$ of side length $2^{-k_1-k_2-3}$ that satisfies $\dist(\D_i,\partial \D_i') \geq 3\cdot 2^{-k_1-k_2-3}$. 

By design of $k_1$, inside each $D_i$ there exists a dyadic cube $S_i$ of side length $2^{-k_1-k_2-3}$ such that $\dist(S_i,\partial D_i) \geq 2^{-k_1-k_2-3}$. Similarly, inside each $\mathcal{S}_{m,i}$ there exists a dyadic cube $S_i'$ of side length $2^{-k_1-k_2-3}$ such that $\dist(S_i',\partial \mathcal{S}_{m,i}) \geq 2^{-k_1-k_2-3}$.

We now construct four bi-Lipschitz homeomorphisms. 

Firstly, for each $i=1,\dots,m$ denote by $D_i'$ the closed $n$-cube which is concentric to $D_i$, contains $D_i$, and $\dist(\partial D_i',D_i) = (3\sqrt{n}d)^{-1}$. In particular, for every $z\in \partial D_i'$, 
\[ (3\sqrt{n}d)^{-1} \leq \dist(z,D_i) \leq (3d)^{-1}.\]
There exists $L_1'$ depending only on $n$ and $d$ such that, for each $i=1,\dots,m$, there exists an $L_1'$-bi-Lipschitz mapping 
\[ f_1 : \left(U\setminus \bigcup_{i=1}^m D_i'\right) \cup \bigcup_{i=1}^m\partial D_i \to \left(U\setminus \bigcup_{i=1}^m D_i'\right) \cup \bigcup_{i=1}^m\partial S_i\] 
with $f_1|_{U\setminus \bigcup_{i=1}^m D_i'} = \text{Id}$, and for each $i=1,\dots,m$, $f_1|_{\partial D_i}$ is an orientation preserving similarity mapping $\partial D_i$ onto $\partial S_i$. Applying Theorem \ref{thm:annulus} on each annulus $\overline{D_i'}\setminus D_i$, $f_1$ can be extended to an $L_1$-bi-Lipschitz mapping 
\[ F_1: \overline{U} \to \overline{\mathcal{U}_0 \setminus \bigcup_{i=1}^m S_i}\] 
with $F_1(\partial D_i) = \partial S_i$, and $L_1$ depending only on $n$ and $d$.

Secondly,  for each $i=1,\dots,m$ let 
\[ \mathcal{S}_{m,i}' = \left[\frac{4i-2m-7/2}{2m+1},\frac{4i-2m-1/2}{2m+1}\right]\times \left[-\frac{3/2}{2m+1},\frac{3/2}{2m+1}\right]^{n-1}\]
so that $\mathcal{S}_{m,i}\subset \mathcal{S}_{m,i}' \subset  (-1,1)^2$ for each $i=1,\dots,m$. Working as above, we can find an $L_2$-bi-Lipschitz map \[F_2 : \mathcal{U}_m \to \overline{\mathcal{U}_0 \setminus \bigcup_{i=1}^n S_i'}\]
 with $F_2(\partial \mathcal{S}_{m,i}) = \partial S_i'$, and $L_2$ depending only on $n$ and $d$.
 
Thirdly, we construct a bi-Lipschitz homeomorphism
\[ F_3 : \overline{\mathcal{U}_0 \setminus \bigcup_{i=1}^n S_i} \to \overline{\mathcal{U}_0 \setminus \bigcup_{i=1}^n \D_i}\]
with $F_3(\partial S_i) = \partial \D_i$.
We describe the steps and leave the technichal details to the reader. Denote by $\mathcal{G}$ the union of all edges of all dyadic cubes of side-length $2^{-k_1-k_2-4}$. Note that the centers of cubes $S_i$ and cubes $\D_i$, $i=1,\dots,m$ are in $\mathcal{G}$.

Let $V_1 = \mathcal{U}_0 \setminus (S_2\cup\cdots\cup S_m)$. There exists a polygonal path $\g_1 \subset \mathcal{G}\cap V_1$ with endpoints the center of $S_1$ and the center of $\D_1$ and 
\[ \dist(\g_1, \partial V_1) \geq 5\cdot 2^{-k_1-k_2-4}.\]
Move the cube $S_1$ inside $V_1$ so that the center of $S_1$ traces $\g_1$ and ends up onto $\D_1$.
Inductively, assume that we have moved the cubes  $S_1,\dots, S_j$ onto the cubes $\D_1,\dots,\D_j$, respectively. Set 
\[ V_{j+1} = \mathcal{U}_0 \setminus (\D_1\cup\cdots\cup \D_j \cup S_{j+2}\cup\cdots\cup S_m).\]
There exists a polygonal path $\g_{j+1} \subset \mathcal{G}\cap V_{j+1}$ with endpoints the center of $S_{j+1}$ and the center of $\D_{j+1}$ and 
\[ \dist(\g_{j+1}, \partial V_{j+1}) \geq 5\cdot 2^{-k_1-k_2-4}.\]
Move the cube $S_{j+1}$ inside $V_{j+1}$ so that the center of $S_{j+1}$ traces $\g_{j+1}$ and ends up onto $\D_{j+1}$. After $m$ steps, we have moved each cube $S_i$ onto the cube $\D_i$ and we obtain an $L_3$-bi-Lipschitz map $F_3$. There are at most $(2^{k_1+k_2+3})!$ different combinations for the position of cubes $S_i$ inside $\mathcal{U}_0$, and there are at most $(2^{k_1+k_2+3})!$ different combinations for the position of cubes $\D_i$ inside $\mathcal{U}_0$. Therefore, $L_3$ depends at most on $n,k_1,k_2$, thus at most on $n,d$.

Fourthly, there exists an $L_4$-bi-Lipschitz homeomorphism 
\[ F_4 : \overline{\mathcal{U}_0 \setminus \bigcup_{i=1}^n S_i'} \to \overline{\mathcal{U}_0 \setminus \bigcup_{i=1}^n \D_i}\]
with $F_3(\partial S_i') = \partial \D_i$ and $L_4$ depending only on $d,\eta$. The construction is similar to that of $F_3$.

The map $F_2^{-1}\circ F_4^{-1}\circ F_3\circ F_1 : \overline{U} \to \overline{\mathcal{U}_0}$ is an $L$-bi-Lipschitz homeomorphism such that $F(\partial D_i) = \partial \mathcal{S}_{m,i}$ for all $i\in\{1,\dots,m\}$, and with $L$ depending only on $n$ and $d$.
\end{proof}
%

We are now ready to prove Proposition \ref{prop:ext-ndim}.

\begin{proof}[{Proof of Proposition \ref{prop:ext-ndim}}]
Since the embedding $f$ can be extended homeomorphically to $\overline{U}$, there exists a domain $U' \subset \R^n$ such that $\partial U' = f(\partial U)$ and $f$ can be extended to a homeomorphism of $\overline{U}$ onto $\overline{U'}$. Since $f$ is a similarity on each boundary component of $U$, every boundary component of $U'$ is the boundary of a cube.

Proposition \ref{prop:ext-ndim} is trivial if $U$ is simply connected. For the rest, we assume that $\partial U$ has at least two components. Moreover, applying two similarities, we may assume that $\diam{U} = \diam{U'}=1$.

Since $f$ is bi-Lipschitz, $U' \in \mathscr{C}_n(d')$ for some $d'$ depending only on $d$ and $L$. By Lemma \ref{lem:standdom}, we may assume that $U=U'=\mathcal{U}_m$, that $f$ maps $\partial \mathcal{U}_0$ onto $\partial \mathcal{U}_0$, and for each $k\in\{1,\dots,m\}$
$f$ maps $\partial \mathcal{S}_{m,k}$ onto $\partial \mathcal{S}_{m,k}$.

If $m=1$, the proposition follows from Theorem \ref{thm:annulus}. Assume now that $m\geq 2$. Let 
\[ \mathcal{S}_0' = \left[\frac{1/2}{2m+1}-1,1-\frac{1/2}{2m+1}\right]^n\] 
and for each $k=1,\dots,m$ let 
\[ \mathcal{S}_{m,k}' = \left[\frac{4k-2m-7/2}{2m+1},\frac{4k-2m-1/2}{2m+1}\right]\times \left[-\frac{3/2}{2m+1},\frac{3/2}{2m+1}\right]^{n-1}\]
so that $\mathcal{S}_{m,k}\subset \mathcal{S}_{m,k}' \subset \mathcal{S}_{0}' \subset (-1,1)^2$ for each $k=1,\dots,m$. Extend $f$ to $\partial \mathcal{S}_0'$ and to each $\mathcal{S}_{m,k}'$ with identity and note that the new embedding, which we still denote by $f$, is $L_1$-bi-Lipschitz with $L_1$ depending only on $L$ and $N$, hence only on $L$, $n$, and $d$.


Applying Theorem \ref{thm:annulus} on the interior of each $\mathcal{S}_{k,m}'\setminus \mathcal{S}_{k,m}$ we obtain $L_1'$-bi-Lipschitz extensions 
$F_k$ on $\overline{\mathcal{S}_{k,m}'\setminus \mathcal{S}_{k,m}}$ with $L_1'$ depending only on $n$, $d$ and $L$. Similarly, we obtain an $L_1''$-bi-Lipschitz extension 
$F_0$ on $[-1,1]^n\setminus \mathcal{S}_{0}'$. The map 
\begin{equation*}
F: \overline{\mathcal{U}_m} \to \overline{\mathcal{U}_m} \qquad\text{with}\quad  F(x) =
\begin{cases}
F_0(x), &\text{ if $x \in \overline{(-1,1)^n \setminus \mathcal{S}_{0}'}$}\\
x, &\text{ if $x \in \overline{\mathcal{S}_{0}' \setminus \bigcup_k\mathcal{S}_{m,k}'}$}\\
F_k(x), &\text{ if $x \in \overline{\mathcal{S}_{m,k}' \setminus \mathcal{S}_{m,k}}$}
\end{cases}
\end{equation*}
is $L'$-BLD for some $L'$ depending only on $n$, $L$ and $d$. By Lemma \ref{lem:BLD}, $f$ is $L'$-bi-Lipschitz.
\end{proof}

\section{Proof of Theorem \ref{thm:main} and Corollary \ref{cor:2}}\label{sec:proof}

Here we prove the following quantitative version of Theorem \ref{thm:main}. The proof of Corollary \ref{cor:2} is given in \textsection\ref{sec:cor}.

\begin{theorem}\label{thm:main2}
Let $n\in\N$ and let $X \subset \R^n$ be a compact $c$-uniformly perfect and $c$-uniformly disconnected set. There exists an $\eta'$-quasisymmetric homeomorphism $F:\R^{n+1} \to \R^{n+1}$ with $F(\mathcal{C}) = X$, and with $\eta'$ depending only on $n$ and $c$.
\end{theorem} 

The main idea of the proof is to divide the set $X$ into appropriate subsets and use the extra dimension of $\R^{n+1}$ to lift each subset to an appropriate height. The subsets and heights have to be chosen carefully to make sure that the lifted subsets are contained in disjoint cubes so that we can appeal to Proposition \ref{prop:ext-ndim}.


We define $\W$ to be the set of finite words formed from the letters $\{1,2\}$, including the empty word $\varepsilon$. Define $\W^N$ to be the set of words in $\W$ whose length is exactly $N$. Given $w\in\mathcal{W}$, we denote by $|w|$ the number of letters that $w$ has, with $|\varepsilon| =0$. 

Let $I_{\varepsilon} = [0,1]$ and given $I_w = [a,b]$ let $I_{w1}=[a,a+\frac13(b-a)]$, $I_{w2} = [b-\frac13(b-a),b]$. For each $w \in \W$, let $\C_w = I_w \cap \C$.

\begin{lemma}\label{lem:cantor}
Let $X$ be a metric space and let $f:\mathcal{C}\to X$ be an $\eta$-quasisymmetric homeomorphism. For each $D>1$, there exists $k\in\N$ depending only on $\eta$ and $D$ with the following property. For any integer $m\geq \log{k}/\log{2}$, there exists a partition $\mathscr{E}_1,\dots,\mathscr{E}_k$ of $ \W^m$
such that 
for any $i\in\{1,\dots,k\}$ and any distinct $w,w' \in \mathscr{E}_i$, 
\begin{equation}\label{eq:cantor7}
\dist(f(\C_w),f(\C_{w'})) \geq D\max\{\diam{f(\C_w)},\diam{f(\C_{w'})}\}.
\end{equation}
\end{lemma}

\begin{proof}
Set $d = (\eta^{-1}((2D)^{-1}))^{-1}$. We show that the lemma holds for $k$ being the integer part of $2d^{\log{2}/\log{3}} +1$. Fix an integer $m\geq \log{k}/\log{2}$. By quasisymmetry of $f$ and (\ref{eq:relQS}), equation (\ref{eq:cantor7}) is satisfied if $\dist(\C_w,\C_{w'}) \geq d3^{-m}$. 

Let $l$ be the integer part of $\log{d}/\log{3} + 1$. By definition, $k \geq 2^l$. For distinct $u,u' \in\W^{m-l}$ we have
\[ \dist(\mathcal{C}_u, \mathcal{C}_{u'}) \geq 3^{-(m-l)} \geq d 3^{-l}.\]
Therefore, for all $w\in \mathcal{W}^m$, there exist at most $2^l$ (hence at most $k$) words $w'\in \W^{m}$ such that $\dist(\mathcal{C}_{w},\mathcal{C}_{w'})\geq d 3^{-l}$.

Let now $\{w_1,\dots,w_{2^m}\}$ be an enumeration of $\W^m$ such that for all $1\leq i<j\leq 2^m$, the set $\C_{w_i}$ lies to the left of the set $\C_{w_j}$. For each $i=1,\dots,k$ define $A_i$ to be the integers in $\{1,\dots,2^m\}$ that are of the form $i+rk$ with $r\in\N\cup\{0\}$ and set $\mathscr{E}_i = \{w_j \colon j\in A_i\}$. It is now straightforward to verify that the sets $\mathscr{E}_i$ satisfy the properties (1) and (2) of the lemma.
\end{proof}

We are now ready to prove Theorem \ref{thm:main2}.

\begin{proof}[{Proof of Theorem \ref{thm:main2}}]
Let $X$ be a compact, $c$-uniformly perfect and $c$-uniformly disconnected subset of $\R^n$. By Theorem \ref{thm:UDQS}, there exists an $\eta$-quasisymmetric homeomorphism $f: \C \to X$ with $\eta$ depending only on $n$ and $c$. 
The first step of the proof is the construction of a bi-Lipschitz mapping $\Phi:\R^{n+1} \to \R^{n+1}$ that unlinks $X$. The second step is the construction of a quasiconformal mapping $G:\R^{n+1} \to \R^{n+1}$ that maps the unlinked image $\Phi(X)$ onto $\C$. The composition $G\circ\Phi$ is the desired map $F$ of Theorem \ref{thm:main}. 

Without loss of generality, we assume that $\diam{X}=1$. For the rest of the proof we write $X_w = f(\mathcal{C}_w)$. Let  $k$ be the number obtained by Lemma \ref{lem:cantor} for $D:=1+8\sqrt{n+1}$ and for $\eta$. Let also 
$N$ be the smallest positive integer such that 
\[3^{-N} \leq \tfrac12\eta^{-1}\left ( (20k\sqrt{n})^{-1}
\right) \qquad\text{and}\qquad N \geq \log{k}/\log{2}.\] 
By (\ref{eq:relQS}) and (\ref{eq:relQS2}), for any two distinct $w,w'\in\W$ with $|w|=|w'|$, and any $u\in\W^N$,
\begin{align}
&(2\eta(3^N))^{-1} \diam{X_w} \leq\diam{X_{wu}} \leq \eta(2\cdot 3^{-N})\diam{X_w}\label{eq:cantor1},\\
&\dist(X_{w},X_{w'}) \geq (2\eta(1))^{-1} \max\{\diam{X_{w}},\diam{X_{w'}}\}.\label{eq:cantor2}
\end{align}

Let $\mathscr{E}^{\varepsilon}_{1}, \dots,\mathscr{E}^{\varepsilon}_k$ be the sets of Lemma \ref{lem:cantor} corresponding to $f$, $D$, and $m=N$. Define $\phi_1: X \to \R$ by 
\[\phi_1|_{X_w}(x) = (4k)^{-1}(i-1),\qquad\text{for $w\in\mathscr{E}^{\varepsilon}_{i}$ and $i\in\{1,\dots,k\}$.}\]

Inductively, suppose that for some $j\in\N$ we have defined $\phi_j : X \to \R$ such that $\phi_j|_{X_{w}}$ is constant whenever $w \in\W^{jN}$. For each $w \in \W^{jN}$, let $\zeta_w : \mathcal{C} \to \mathcal{C}_w$ be a similarity and note that $f|_{\mathcal{C}_{w}} \circ \zeta_w : \mathcal{C} \to X_w$ is $\eta$-quasisymmetric. Let $\mathscr{E}^{w}_1, \dots, \mathscr{E}^{w}_k$ be the sets of $\W^N$ from Lemma \ref{lem:cantor} applied to $f|_{\mathcal{C}_{w}} \circ \zeta_w$, $D$, and $m=N$. 
Define $\phi_{j+1}: X \to \R$ such that
\[ \phi_{j+1}|_{X_{wu}} (x) = \phi_j|_{X_{w}}(x) + (4k)^{-1}(i-1)\diam{X_{w}},\]
where $w\in\W^{jN}$, $u\in \mathscr{E}^{w}_i$ and $i\in \{1,\dots,k\}$.

By (\ref{eq:cantor1}), for positive integers $j< j'$ we have
\begin{align*} 
\|\phi_j-\phi_{j'}\|_{\infty} \leq \sum_{l=j}^{j'-1} \|\phi_l-\phi_{l+1}\|_{\infty} &\leq \sum_{l=j}^{j'-1}(4k)^{-1}(k-1)\diam{X_w}\\ 
&\leq \sum_{l=j}^{j'-1}\tfrac14\left(\eta(2\cdot 3^{-N})\right)^{l}\\
&\leq \frac{\left(\eta(2\cdot 3^{-N})\right)^{j}}{4-4\eta(2\cdot 3^{-N})}
\end{align*}
which, by choice of $N$, goes to $0$ as $j\to\infty$. Therefore, the mappings $\phi_j$ converge uniformly to a mapping $\phi: X \to \R$.

We claim that $\phi$ is Lipschitz with the Lipschitz constant depending only on $n$ and $c$. To see that, let $x,y \in X$ and let $j\in\N$ be the unique integer such that there exists $w \in \W^{jN}$ and there exist distinct $u,u' \in\W^N$ with $x \in X_{wu}$ and $y\in X_{wu'}$. On the one hand, by (\ref{eq:cantor1}), (\ref{eq:cantor2}),
\begin{align*}
|x-y| \geq \dist(X_{wu},X_{wu'}) \geq \frac1{2\eta(1)}\min_{v\in\W^N}\diam{X_{wv}} \geq \frac{\diam{X_{w}}}{4\eta(1)\eta(3^N)}.
\end{align*}
On the other hand, by (\ref{eq:cantor1}) and the fact that $\phi_{j}(x)=\phi_{j}(y)$,
\begin{align*}
|\phi(x)-\phi(y)| &\leq |\phi(x)-\phi_{j}(x)| + |\phi_{j}(x)-\phi_{j}(y)| + |\phi(y)-\phi_{j}(y)|\\
&\leq 2\sum_{l=j}^{\infty}\|\phi_l-\phi_{l+1}\|_{\infty}\\
&\leq 2(4k)^{-1}(k-1)\diam{X_{w}}\sum_{l=0}^{\infty}(\eta(2\cdot 3^{-N}))^l\\
&\leq \frac{1}{2-2\eta(2\cdot 3^{-N})}\diam{X_{w}}.
\end{align*}
Therefore,
\[ |\phi(x)-\phi(y)| \leq \frac{2\eta(1)\eta(3^N)}{1-\eta(2\cdot 3^{-N})} |x-y|\]
and the claim follows,

Fix $x_0\in X$, $B_0 = B^n(x_0,5\diam{X})$ and set $\phi|_{\R^n\setminus B_0} \equiv 0$. Then, the map 
\[ \phi: (\R^n \setminus B_0) \cup X \to \R\] 
is $L$-Lipschitz for some $L$ depending only on $n,c$ and, by Kirszbraun Theorem, there exists an $L$-Lipschitz extension of $\phi$ to $\R^n$ which we also denote by $\phi$. Then, the mapping 
\[ \Phi:\R^{n+1} \to \R^{n+1}, \qquad \text{defined by $\Phi(x,z) = (x,\phi(x)+z)$}\] 
is $L'$-bi-Lipschitz with $L'=2(L+2)$ \cite[Lemma 5.3.2]{Vthesis}.

For each $j=0,1,\dots$ and each $w\in \W^{jN}$ fix $x_w \in X_{w}$ and set 
\[\mathsf{K}_w := x_w + [-2\diam{X_w}, 2\diam{X_w}]^n = \{x_w+(2\diam{X_w})y : y \in  [-1, 1]^n\}.\]
We claim that if $w \in \W^{jN}$ and $u\in \W^N$, then 
\begin{equation}\label{eq:cantor8}
\mathsf{K}_{wu} \subset \mathsf{K}_{w} \qquad\text{and}\qquad \dist(\mathsf{K}_{wu}, \partial \mathsf{K}_{w}) \geq \frac12\diam{X_w}.
\end{equation} 
Indeed, if $x\in \mathsf{K}_{wu}$, then by (\ref{eq:cantor1}) and the choice of $N$,
\begin{align*}
|x-x_w| \leq |x-x_{wu}| + |x_{wu}-x_w| &\leq 2\sqrt{n}\diam{X_{wu}} + \diam{X_w}\\ 
&\leq (1+2\sqrt{n}\eta(2\cdot 3^{-N}))\diam{X_w}\\
&< \frac32\diam{X_w}
\end{align*}
which proves both claims.

We remark that if $w,w' \in \W^{jN}$ are distinct, then $\mathsf{K}_{w}$ may intersect $\mathsf{K}_{w'}$. This is why we lift different sets $X_w$ to different heights. 

For each $j=0,1,\dots$ and each $w\in\W^{jN}$ define 
\[ \mathcal{K}_w := \mathsf{K}_w\times[\phi_j(x_w)-2\diam{X_w},\phi_j(x_w)+2\diam{X_w}] \subset \R^{n+1}. \]

We first claim that for all $j\in\N$, for all $w \in \W^{jN}$ and for all $u\in \W^{N}$,
\begin{equation}\label{eq:cantor5}
\mathcal{K}_{wu} \subset \mathcal{K}_{w} \qquad\text{and}\qquad \dist(\mathcal{K}_{wu}, \partial \mathcal{K}_{w}) \geq (8\sqrt{n+1})^{-1}\diam{\mathcal{K}_w}.
\end{equation}
By (\ref{eq:cantor8}), it suffices to check only the $(n+1)$-th coordinate. Let $z = (z_1,\dots,z_{n+1})\in \mathcal{K}_{wu}$. By (\ref{eq:cantor1}) and the choice of $N$, 
\begin{align*}
|z_{n+1} - \phi_j(x_w)| &\leq |z_{n+1}-\phi_{j+1}(x_{wu})| + |\phi_j(x_w)-\phi_{j+1}(x_{wu})|\\ 
&\leq 2\diam{X_{wu}} + (4k)^{-1}(k-1)\diam{X_w}\\ 
&\leq  2\eta(2\cdot 3^{-N})\diam{X_w} + \tfrac14\diam{X_w}\\
&\leq \frac12\diam{X_w}
\end{align*} 
and the claim follows from (\ref{eq:cantor7}).

Secondly, we claim that for all integers $j\geq 0$, all $w\in \mathcal{W}^{jN}$ and all distinct $u, u' \in \W^{N}$,
\begin{equation}\label{eq:cantor4}
\dist(\mathcal{K}_{wu},\mathcal{K}_{wu'}) \geq (8\eta(3^N)\sqrt{n+1})^{-1}\diam{\mathcal{K}_w}.
\end{equation}
To prove (\ref{eq:cantor4}), 
let $x\in\mathcal{K}_{wu}$ and $x'\in\mathcal{K}_{wu'}$. There are two cases to consider.

\emph{Case 1.} Suppose that $u,u' \in \mathscr{E}_i^{w}$. Then $\phi_{j+1}(x_{wu}) = \phi_{j+1}(x_{wu'})$ and by the choice of $D$ we have
\begin{align*}
|x-x'| &\geq |x_{wu} - x_{wu'}| - \diam{\mathcal{K}_{wu}} - \diam{\mathcal{K}_{wu'}}\\
&\geq \dist(X_{wu},X_{wu'}) - 2\max\{\diam{\mathcal{K}_{wu}},\diam{\mathcal{K}_{wu'}}\}\\
&\geq (D-8\sqrt{n+1})\max\{\diam{X_{wu}},\diam{X_{wu'}}\}\\
&=\max\{\diam{X_{wu}},\diam{X_{wu'}}\}.
\end{align*}

\emph{Case 2.} Suppose that $u \in \mathscr{E}_i^w$ and $u' \in \mathscr{E}_{i'}^w$ with $i\neq i'$. By the choice of $N$, we calculate the vertical difference of $x,x'$
\begin{align*}
|x-x'| &\geq |\phi_{j+1}(x_{w_0u}) - \phi_{j+1}(x_{wu'})| - 2\diam{X_{wu}} - 2\diam{X_{wu}}\\
&\geq (4k)^{-1}\diam{X_{w}} - 4\max\{\diam{X_{wu}},\diam{X_{wu'}}\}\\
&\geq \left( \frac{1}{4k\eta(2\cdot 3^{-N})} - 4\right)\max\{\diam{X_{wu}},\diam{X_{wu'}}\}\\
&\geq \max\{\diam{X_{wu}},\diam{X_{wu'}}\}.
\end{align*}
In either case, (\ref{eq:cantor4}) follows now from (\ref{eq:cantor1}).

Thirdly, by (\ref{eq:cantor1}), we have that for all $j\in\N$, for all $w \in \W^{jN}$ and for all $u\in \W^{N}$,
\begin{equation}\label{eq:cantor7}
\diam{\mathcal{K}_{wu}} \geq (2\eta(3^N))^{-1}\diam{\mathcal{K}_{w}}.
\end{equation}

Finally, by design of $\Phi$,
\begin{equation}\label{eq:cantor6}
\mathcal{K}_{w} \cap \Phi(X) = \Phi(X_w) \quad\text{and}\quad \dist(\Phi(X_w),\partial\mathcal{K}_w) \geq \frac12\diam{X_w}.
\end{equation}

For each $j=0,1,\dots$ and $w\in\W^{jN}$, let $z_w$ be the centre of $I_w$ and define the cube
\[ \mathcal{Q}_w = \left [z_w-\tfrac56 3^{-jN}, z_w+\tfrac563^{-jN}\right]\times \left[-\tfrac56 3^{-jN}, \tfrac563^{-jN}\right]^n.\]
For each $w\in\W^{jN}$, let $g_w :\partial\mathcal{K}_{w}\to\partial\mathcal{Q}_{w}$ be an orientation preserving similarity map. By (\ref{eq:cantor5}), (\ref{eq:cantor4}) and (\ref{eq:cantor7}), Proposition \ref{prop:ext-ndim} applies and there exists $\Lambda>1$ depending only on $n,\eta$, (hence only on $n,c$) and there exists $G:\R^{n+1} \to \R^{n+1}$ such that 
\begin{enumerate}
\item $G$ is the identity outside of $B_0$ and
\item for all $w\in\W^{jN}$, the restriction of $G$ on $\mathcal{K}_{w}\setminus (\bigcup_{u\in\W^N}\mathcal{K}_{wu})$ extends $g_w$ and is a $(\frac{\diam{\mathcal{Q}_w}}{\diam{\mathcal{K}_w}},\Lambda)$-quasisimilarity that maps $\mathcal{K}_{w}\setminus (\bigcup_{u\in\W^N}\mathcal{K}_{wu})$ onto $\mathcal{Q}_{w}\setminus (\bigcup_{u\in\W^N}\mathcal{Q}_{wu})$.
\end{enumerate}

Therefore, by a theorem of V\"ais\"al\"a on removability of singularities \cite[Theorem 35.1]{Vais1}, $G$ is $K$-quasiconformal with $K$ depending only on $c$ and $n$. Set $F = G\circ\Phi$ and note that $F$ extends $f$. Therefore, $F(X) = \C$.
\end{proof}

\subsection{Proof of Corollary \ref{cor:2}}\label{sec:cor}

For the proof of Corollary \ref{cor:2}, recall that a set $E \subset [0,1]$ is \emph{porous} if there exists $c\geq 1$ such that for any interval $I\subset [0,1]$, there exists an interval $J\subset I\setminus X'$ of length $|J|\geq c^{-1}|I|$.

\begin{proof}[{Proof of Corollary \ref{cor:2}}]
Let $X\subset \R^n$ be bounded and $c$-uniformly disconnected. We show that $X$ is contained in a compact uniformly perfect and uniformly disconnected set. Replacing $X$ by $\overline{X}$, we may assume that $X$ is compact. By Theorem 1 in \cite{MM2} (see also \cite[p. 275]{MM2} for discussion and \cite[Theorem 3.8]{BV} for a more general statement), there exists an $\eta$-quasisymmetric map $f:[0,1] \to \S^n$ such that $X\subset [0,1]$, $\{0,1\}\subset X$, and $\eta$ depending only on $c$. Therefore, we may assume that $X \subset [0,1]$. For each component $I$ of $[0,1]\setminus X$, let $\phi_I : \R \to \R$ be a similarity that maps $[0,1]$ onto $\overline{I}$. Define now
\[ X' := X\cup \bigcup_I\phi_I(\mathcal{C}) = \overline{ \bigcup_I\phi_I(\mathcal{C}) }.\]
where the union is over all components $I$ of $[0,1]\setminus X$. It remains to show that $X'$ is uniformly perfect and uniformly disconnected. Recall that $\mathcal{C}$ is $C_0$-uniformly perfect and $C_0$-uniformly disconnected for some $C_0> 1$. 

To show that $X'$ is uniformly perfect, fix $x\in X'$ and $r\in (0,1)$. We claim that there exists universal $C\geq 1$ such that 
\[ X'\cap \left((x-r,x+r)\setminus (x-r/C,x+r/C) \right)\neq \emptyset.\]
If  $(x-r,x+r) \cap X = \emptyset$, then $(x-r,x+r) \cap X' = (x-r,x+r) \cap \mathcal{C}$ and the claim is true for $C=C_0$.  Suppose now that $(x-r,x+r) \cap X \neq \emptyset$. If there exists $z\in (x-r,x+r) \cap X $ with $|z-x| \geq r/2$, then the claim is true for $C=2$. Suppose now that $(x-r/2,x+r/2) \cap X = (x-r,x+r) \cap X \neq \emptyset$ and let $z \in(x-r/2,x+r/2) \cap X $ such that $|z-x|$ is maximal. Since $z$ is the endpoint of some component $I$ of $[0,1]\setminus X$, there exists a rescaled copy $Z$ of $\mathcal{C}$ such that $z\in Z$, $Z\subset X'$ and $\diam{Z}\geq r/6$. Then the claim holds true for $C = 6$.

To show that $X'$ is uniformly disconnected, by \cite[Theorem 1]{MM2}, we need to show that $X'$ is porous. Fix an interval $I\subset [0,1]$. The porosity of $X$ implies that there exists an interval $J' \subset I\setminus X$ such that $|J'|\geq |I|/c_1$ where $c_1$ depends only on $C$. Now, the porosity of $J'\cap X'$ (since it is a subset of a copy of $\mathcal{C}$) implies that there exists an interval $J \subset J'\setminus X'$ such that $|J|\geq |J'|/c_0$ where $c_0$ depends only on $C_0$. Alltogether, $J\subset I\setminus X'$ and $|J|\geq (c_1c_0)^{-1}|I|$.
\end{proof}

\bibliographystyle{alpha}

\bibliography{proposal}
\end{document}